\documentclass[11pt,fullpage]{article}

\newcommand{\ol}{\setlength{\itemsep}{0pt.}\begin{enumerate}}
\newcommand{\eol}{\end{enumerate}\setlength{\itemsep}{-\parsep}}
\newcommand{\ignore}[1]{}
\begin{document}
\begin{center}
{\bf  Combinatorial and Algorithmic Aspects of Hyperbolic Polynomials  
 }

\vskip 4pt
{Leonid Gurvits }\\
\vskip 4pt
{\tt gurvits@lanl.gov}\\
\vskip 4pt
Los Alamos National Laboratory, Los Alamos , NM 87545 , USA.
\end{center}




\begin{abstract}
Let $p(x_1,...,x_n) =\sum_{ (r_1,...,r_n) \in I_{n,n} } a_{(r_1,...,r_n) } \prod_{1 \leq i \leq n} x_{i}^{r_{i}}$ 
be homogeneous polynomial of degree $n$ in $n$ real variables with integer nonnegative coefficients. 
The support of such polynomial $p(x_1,...,x_n)$  
is defined as $supp(p) = \{(r_1,...,r_n) \in I_{n,n} : a_{(r_1,...,r_n)} \neq 0 \}$ . The convex hull $CO(supp(p))$ of $supp(p)$ is called
the Newton polytope of $p$ .
We study the following decision problems , which are far-reaching generalizations of the classical perfect matching problem  :
\begin{itemize}
\item
{\bf Problem 1 .} Consider a homogeneous polynomial $p(x_1,...,x_n)$ of degree $n$ in $n$ real variables with
nonnegative integer coefficients given as a black box (oracle ) .
{\it Is it true that $(1,1,..,1) \in supp(p)$ ? }
\item
{\bf Problem 2 .} Consider a homogeneous polynomial $p(x_1,...,x_n)$ of degree $n$ in $n$ real variables with
nonnegative integer coefficients given as a black box (oracle ) .
{\it Is it true that $(1,1,..,1) \in CO(supp(p))$ ? }
\end{itemize}
We prove that for hyperbolic polynomials these two problems are equivalent and can be solved by deterministic polynomial-time oracle algorithms .
This result is based on a "hyperbolic" generalization of the Rado theorem  .
We also present combinatorial and algebraic applications of this "hyperbolic" generalization of the Rado theorem 
(prove that the support $supp(p)$ of $P$-hyperbolic polynomial $p$ is an intersection of some Integral Polymatroid with the hyperplane
$\{(r_1,...,r_n) : \sum_{1 \leq i \leq n} r_i = n \}$ )
and pose some open problems.

\end{abstract} 
 
 \newtheorem{ALGORITHM}{Algorithm}[section]
\newenvironment{algorithm}{\begin{ALGORITHM} \hspace{-.85em} {\bf :} 
}%
                        {\end{ALGORITHM}}
\newtheorem{THEOREM}{Theorem}[section]
\newenvironment{theorem}{\begin{THEOREM} \hspace{-.85em} {\bf :} 
}%
                        {\end{THEOREM}}
\newtheorem{LEMMA}[THEOREM]{Lemma}
\newenvironment{lemma}{\begin{LEMMA} \hspace{-.85em} {\bf :} }%
                      {\end{LEMMA}}
\newtheorem{COROLLARY}[THEOREM]{Corollary}
\newenvironment{corollary}{\begin{COROLLARY} \hspace{-.85em} {\bf 
:} }%
                          {\end{COROLLARY}}
\newtheorem{PROPOSITION}[THEOREM]{Proposition}
\newenvironment{proposition}{\begin{PROPOSITION} \hspace{-.85em} 
{\bf :} }%
                            {\end{PROPOSITION}}
\newtheorem{DEFINITION}[THEOREM]{Definition}
\newenvironment{definition}{\begin{DEFINITION} \hspace{-.85em} {\bf 
:} \rm}%
                            {\end{DEFINITION}}
\newtheorem{EXAMPLE}[THEOREM]{Example}
\newenvironment{example}{\begin{EXAMPLE} \hspace{-.85em} {\bf :} 
\rm}%
                            {\end{EXAMPLE}}
\newtheorem{CONJECTURE}[THEOREM]{Conjecture}
\newenvironment{conjecture}{\begin{CONJECTURE} \hspace{-.85em} 
{\bf :} \rm}%
                            {\end{CONJECTURE}}
\newtheorem{PROBLEM}[THEOREM]{Problem}
\newenvironment{problem}{\begin{PROBLEM} \hspace{-.85em} {\bf :} 
\rm}%
                            {\end{PROBLEM}}
\newtheorem{QUESTION}[THEOREM]{Question}
\newenvironment{question}{\begin{QUESTION} \hspace{-.85em} {\bf :} 
\rm}%
                            {\end{QUESTION}}
\newtheorem{FACT}[THEOREM]{Fact}
\newenvironment{fact}{\begin{FACT} \hspace{-.85em} {\bf :} 
\rm}%
                            {\end{FACT}}

\newtheorem{REMARK}[THEOREM]{Remark}
\newenvironment{remark}{\begin{REMARK} \hspace{-.85em} {\bf :} 
\rm}%
                            {\end{REMARK}}
\newcommand{\alg}{\begin{algorithm}} 
\newcommand{\thm}{\begin{theorem}}
\newcommand{\lem}{\begin{lemma}}
\newcommand{\pro}{\begin{proposition}}
\newcommand{\dfn}{\begin{definition}}
\newcommand{\fac}{\begin{fact}}

\newcommand{\rem}{\begin{remark}}
\newcommand{\xam}{\begin{example}}
\newcommand{\cnj}{\begin{conjecture}}
\newcommand{\prb}{\begin{problem}}
\newcommand{\que}{\begin{question}}
\newcommand{\cor}{\begin{corollary}}
\newcommand{\prf}{\noindent{\bf Proof:} }
\newcommand{\ethm}{\end{theorem}}
\newcommand{\elem}{\end{lemma}}
\newcommand{\epro}{\end{proposition}}
\newcommand{\edfn}{\bbox\end{definition}}
\newcommand{\efac}{\bbox\end{fact}}

\newcommand{\erem}{\bbox\end{remark}}
\newcommand{\exam}{\bbox\end{example}}
\newcommand{\ealg}{\end{algorithm}}
\newcommand{\ecnj}{\bbox\end{conjecture}}
\newcommand{\eprb}{\bbox\end{problem}}
\newcommand{\eque}{\bbox\end{question}}
\newcommand{\ecor}{\end{corollary}}
\newcommand{\eprf}{\bbox}
\newcommand{\beqn}{\begin{equation}}
\newcommand{\eeqn}{\end{equation}}
\newcommand{\wbox}{\mbox{$\sqcap$\llap{$\sqcup$}}}
\newcommand{\bbox}{\vrule height7pt width4pt depth1pt}
\newcommand{\qed}{\bbox}
\def\sup{^}
\def\Tp{Tchebyshef polynomial}
\def\Tps{TchebysDeto be the maximafine $A(n,d)$ l size of a code with distance 
$d$hef polynomials}
\newcommand{\rarrow}{\rightarrow}
\newcommand{\larrow}{\leftarrow}
\newcommand{\grad}{\bigtriangledown}

\overfullrule=0pt
\def\setof#1{\lbrace #1 \rbrace}

\section{Introduction and motivating examples}
{ \bf The layout of the paper :} \\
{\it We introduce the main topics and motivations in Section 1 . In Section 1.1
we present a naive algorithm to solve Problem 1 in the general case .
We show in Appendix {\bf D} that this algorithm is , in a sense , optimal . \\
(Incidentally (or not) , the situation here is very similar with
the optimality of the square root in the famous quantum Grover's search algorithm .)

\noindent In Section 1.2 we remind the basic properties of hyperbolic polynomials
used in this paper .

\noindent In Section 2 we state a hyperbolic analogue of the Rado theorem (Theorem 2.2) , which is the main mathematical
result of the paper . Theorem 2.2 sheds more light on the algebraic-geometric nature
of such fundamental combinatorial results as Hall's and Rado's theorems .\\ 
In Section 2.1 we define and study doubly-stochastic polynomials 
(an useful generalization of standard doubly-stochastic matrices ). We also state there
a hyperbolic analogue of the van der Waerden conjecture .\\

\noindent In Section 3 we introduce and analyse the ellipsoid algorithm which solves Problem 1 and Problem 2
on the class of $S$-hyperbolic polynomials . The essence of the results in Section 3 is that
once Hall's like conditions (the exponential number of them) are proved to be necessary and sufficient ,
they can be checked by a polynomial time deterministic oracle algorithms . The algorithm , which
we use , is based not on the linear programming but on some rather
nonlinear convex programs similar to considered in \cite{nr} , \cite{GS} , \cite{GS1}.\\

In section 4 we introduce and analyse another algorithm ,
which is a "polynomial" generalization of the Sinkhorn Scaling .\\

\noindent In Section 5 we use Theorem 2.2 to get a more refine (polymatroidal) properties of the supports
of $P$-hyperbolic polynomials and explain how our results generalize the main result from \cite{half}.\\

\noindent In Section 6 we pose some open problems and share our enthusiasm about the topic .

\noindent The proofs of the main results are presented in Appendices {\bf A,B,C,D} . }

Let $p(x_1,...,x_n) =\sum_{ (r_1,...,r_n) \in I_{n,n} } a_{(r_1,...,r_n) } \prod_{1 \leq i \leq n} x_{i}^{r_{i}}$ 
be homogeneous polynomial of degree $n$ in $n$ real variables. Here $I_{k,n}$ stands for the set of vectors
$r = (r_1,...,r_k)$ with nonnegative integer components and  $\sum_{1 \leq i \leq k} r_i = n$.
In this paper we primarily study homogeneous polynomials with nonnegative integer  coefficients .
\dfn
The support of the polynomial $p(x_1,...,x_n)$ as above 
is defined as $supp(p) = \{(r_1,...,r_n) \in I_{n,n} : a_{(r_1,...,r_n)} \neq 0 \}$ . The convex hull $CO(supp(p))$ of $supp(p)$ is called
the Newton polytope of $p$ .
\edfn
We will study the following decision problems :
\begin{itemize}
\item
{\bf Problem 1 .} Consider a homogeneous polynomial $p(x_1,...,x_n)$ of degree $n$ in $n$ real variables with
nonnegative integer coefficients given as a black box (oracle ) .
{\it Is it true that $(1,1,..,1) \in supp(p)$ ? }\\
An equivalent question is : {\it Is it true that $\frac{\partial^n}{\partial x_1...\partial x_n} p(x_1,...,x_n) \neq 0$ ? }\\
\item
{\bf Problem 2 .} Consider a homogeneous polynomial $p(x_1,...,x_n)$ of degree $n$ in $n$ real variables with
nonnegative integer coefficients given as a black box (oracle ) .
{\it Is it true that $(1,1,..,1) \in CO(supp(p))$ ? }
\end{itemize}
Our goal is solve these decision problems using deterministic polynomial-time oracle algorithms , i.e.
algorithms which evaluate the given polynomial $p(.)$ at a number of rational vectors $(q_1,...,q_n)$ which is 
polynomial in $n$ and $\log(p(1,1,..,1))$; these rational vectors $(q_1,...,q_n)$ are supposed to 
have bit-wise complexity which is polynomial in $n$ and $\log(p(1,1,..,1))$ ; and the additional auxilary arithmetic 
computations also take a polynomial number of steps in $n$ and $\log(p(1,1,..,1))$ .\\ 

The next example explains some (well known ) origins of the both problems .
\xam
Consider first the following homogeneous polynomial from \cite{Bar4} : $p(x_1,...,x_n) = tr (( D(x) A)^{n})$ ,
where $D(x)$ is a $n \times n$ diagonal matrix $Diag(x_1,...,x_n)$ ; and $A$ is $n \times n$ matrix
with $(0,1)$ entries , i.e. $A$ is an adjacency matrix of some directed graph $\Gamma$ . Clearly , this polynomial
$p(x_1,...,x_n)$ has nonnegative integer coefficients . It was proved in  \cite{Bar4} that
$ \frac{1}{n} \frac{\partial^n}{\partial x_1...\partial x_n} tr ( (D(x) A)^{n})$ is equal to the number of Hamiltonian
circuits in the graph $\Gamma$ . Notice that the polynomial $tr ( D(x) A)^{n}$ can be evaluated in $O(n^3 \log(n))$
arithmetic operations and $(1,1,...,1) \in supp(p)$ iff there exists a Hamiltonian circuit in the graph $\Gamma$.
Also $\log(p(1,1,..,1)) \leq n \log(n)$.
Therefore , unless $P=NP$ , there is no hope to design  deterministic polynomial-time oracle algorithm solving Problem 1 in
this case . (The author is indebted to A.Barvinok for pointing out this polynomial . ) \\
Consider , with the same adjacency matrix $A$ , another homogeneous polynomial \\
$Mul(x_1,...,x_n) = \prod_{1 \leq i \leq n} \sum_{1 \leq i \leq n} A(i,j) x_j$ .
Then  \\
$\frac{\partial^n}{\partial x_1...\partial x_n} Mul(x_1,...,x_n) = Per(A)$ .\\
Therefore for the multilinear polynomial $Mul(x_1,...,x_n)$ Problem 1 is a "black box"
analogue of checking the existence of the perfect bipartite matching . \\

Next consider the following class of determinantal polynomials  :
$$
q(x_1,...,x_n) = \det( \sum _{1 \leq i \leq n} A_i x_i) ,
$$ 
where ${\bf A} = (A_1,...,A_n)$ is a $n$-tuple of positive semidefinite $n \times n$ hermitian matrices , i.e. $A_i \succeq 0$ ,  with integer entries .
Recall that the mixed discriminant 
$$
D({\bf A}) = \frac{\partial^n}{\partial \alpha_1...\partial \alpha_n} \det( \sum _{1 \leq i \leq n} A_i x_i) .
$$
(If the matrices $A_i$ above are diagonal , i.e. $A_i = Diag(b_{i,1},...,b_{i,n}) , 1 \leq i \leq n$ ,\\
then the mixed discriminant is reduced to the permanent : $D({\bf A}) = Per(B) , B = \{b_{i,j} ; 1 \leq i,j \leq n$ ).\\ 

It is well known (see , for instance ,  \cite{GS1} ) that a determinantal polynomial $q(.)$ can be represented as
\beqn
q(x_1,..,x_n) = \sum_{r \in I_{n,n}} \prod_{1 \leq i \leq n } x_{i}^{r_{i}} D({\bf A}_{r}) 
\frac{1}{\prod_{1 \leq i \leq n } r_{i}!} ,
\eeqn
where a $n$-tuple ${\bf A}_{r}$ of square matrices consists of $r_i$ copies of $A_{i}, 1 \leq i \leq k $ .
One of the equivalent formulations \cite{panov} of the classical Rado theorem states that $D({\bf A}_{(1,1,...,1)}) > 0$ iff
\beqn
Rank(\sum_{i \in S} A_i) \geq |S| \mbox{\ for all } S \subset \{1,2,...,n \} 
\eeqn
(The diagonal case is the famous Hall's theorem on the perfect bipartite matchings .) \\
The Rado theorem is just a particular case of famous Edmonds theorem on the rank of intersection of two matroids .
Therefore , given a $n$-tuple ${\bf A} = (A_1,...,A_n)$ of positive semidefinite $n \times n$ hermitian matrices , one can decide
in deterministic polynomial time if $D({\bf A}) > 0$ . We will explain below that this decision problem can be solved by a
deterministic polynomial-time oracle algorithm . I.e. we only use some values of $\det( \sum _{1 \leq i \leq n} A_i x_i)$ without
reconstructing the actual tuple ${\bf A} = (A_1,...,A_n)$ . \\

{\bf The natural question , in our opinion , is which algebraic-geometric properties make the class of determinantal
polynomials "easy" and the class of Barvinok's polynomials $tr ( D(x) A)^{n}$ "hard" . This paper suggests one answer to the question. } \\

One important corollary of the Rado conditions (2) is that
\beqn
supp(q) = CO(supp(q)) \cap I_{n,n}.
\eeqn
I.e. if integer vectors $r,r(1),r(2),...,r(k) \in I(n,n)$ and 
$$ r = \sum_{1 \leq i \leq k } a(i)r(i) ,
a(i) \geq 0 , 1 \leq i \leq k ; \sum_{1 \leq i \leq k } a(i) ,
$$
and $D({\bf A}_{r(i)}) > 0 , 1 \leq i \leq k $ then also $D({\bf A}_{r}) > 0$ .
Notice that in this case Problem 1 and Problem 2 are equivalent . \\
We can rewrite  Rado conditions (2) as follows :
\beqn
\max_{r \in supp(q)}\sum_{i \in S} r_i \geq |S| \mbox{\ for all } S \subset \{1,2,...,n \}
\eeqn

Putting things together we get the following Fact . 
\fac
The following properties of determinantal polynomial
$q(x_1,...,x_n) = \det( \sum _{1 \leq i \leq n} A_i x_i) $ with $n \times n$ hermitian matrices $A_i \succeq 0 ,1 \leq i \leq n $ are equivalent .
\begin{enumerate}
\item
$(1,1,..,1) \notin supp(q)$.
\item
$(1,1,..,1) \notin CO(supp(q))$.
\item
There exists nonempty $S \subset \{1,2,...,n \}$ such that
\beqn
\sum _{1 \leq i \leq n} r_i s_i < \sum _{1 \leq i \leq n}  s_i = |S| \mbox{\ for all} (r_1,...,r_n) \in supp(q) ,
\eeqn 

where $(s_1,...,s_n)$ is a characteristic function of the subset $S$ , i.e. $s_i = 1$ if $i \in S$ , and $s_i = 0$ otherwise . \\

{\it Notice that if (5) holds then the distance $dist(e,CO(supp(q)))$  from the vector $e = (1,...,1)$ to the Newton polytope \\ 
$CO(supp(q))$ is at least $\sqrt{\frac{n}{|S| (n- |S|) }} \geq \frac{2}{\sqrt{n}}$ } .
\end{enumerate}
\efac 
We will show that for any class of polynomials satisfying Fact 1.3 there exists a deterministic polynomial-time oracle algorithm
solving both Problem 1 and Problem 2 , which are , of course , equivalent in this case . Our algorithm is based on
the reduction to some convex programming problem and the consequent use of the Ellipsoids method . \\
The next fact about determinantal polynomials , namely their hyperbolicity , is "responsible" for Fact 1.3 .

\fac
Consider a determinantal polynomial $q((x_1,...,x_n) = \det( \sum _{1 \leq i \leq n} A_i x_i) $ with $A_i \succeq 0 ,1 \leq i \leq n $ .
Assume that $q$ is not identically zero , i.e. that $B = :\sum _{1 \leq i \leq n} A_i  \succ 0 $ (the sum is strictly positive definite ).
For a real vector $(x_1,...,x_n) \in R^n$ consider the following polynomial equation of degree $n$ in one variable :
\beqn
P(t) = q(x_1 - t ,x_2 - t , ..., x_n - t) = \det( \sum _{1 \leq i \leq n} A_i x_i - t \sum _{1 \leq i \leq n} A_i ) = 0 .
\eeqn
Equation (6) has $n$ real roots roots counting the multiplicities ; if the real vector $(x_1,...,x_n) \in R^n$ has
nonnegative entries then all roots of (6) are nonnegative real numbers .
\efac

The main result of this paper that this hyperbolicity , which we will describe formally in Section 1.1 , is
sufficient for Fact 1.3 ; i.e. Fact 1.4 implies Fact 1.3 .
\exam
\subsection{"Naive" algorithms}
One possible "naive" algorithm to solve Problem 1 is just to compute $\frac{\partial^n}{\partial x_1...\partial x_n} p(x_1,...,x_n)$ .
Recall that the number of coefficients of a homogeneous polynomial of degree $n$ in $n$ real variables is equal to $\frac{(2n-1)!}{n! (n-1)!} 
\approx 2^{2n}$ . We can compute all the coefficients of $p(x_1,...,x_n)$ via standard multidimensional interpolation ,
but this interpolation will need $\frac{(2n-1)!}{n! (n-1)!}\approx 2^{2n}$ oracle calls . There is an algorithm which computes
$\frac{\partial^n}{\partial x_1...\partial x_n} p(x_1,...,x_n)$ using only $2^{n-1}$ oracle calls :
\begin{eqnarray} 
&\frac{\partial^n}{\partial x_1... \partial x_N} p(x_1,...,x_n)=  \\ \nonumber
&2^{-n+1} \sum_{b_{i} \in \{-1, +1 \}, 2 \leq i \leq n } 
 p(1,b_{2},...,b_{n}) \prod_{2 \leq i \leq n } b_i .
\end{eqnarray}
This formula is , in a sense , optimal . I.e. there exists a nearly matching lower bound .
The corresponding result and connections to computations of the permanent are presented in Appendix D . \\

We will explain below that if $p$ is a homogeneous polynomial  of degree $n$ in $n$ real variables with
nonnegative integer coefficients then
$(1,1,..,1) \in CO(supp(p))$ iff $p(x_1,...,x_n) \geq \prod_{1 \leq i \leq n} x_i$ for
all vectors $(x_1,...,x_n)$ with positive real coordinates .
Therefore Problem 2 is equivalent to checking if the polynomial
$P(y_1,...,y_n) = p(1+y_1^{2} ,...,1+y_n^{2}) -  \prod_{1 \leq i \leq n} 1+y_i^{2}$ is nonnegative on $R^{n}$.

\subsection{Hyperbolic polynomials}
The following concept  of hyperbolic polynomials was originated in the theory of partial differential equations \cite{gar}, \cite{horm} ,\cite{kry} .\\
A homogeneous polynomial $p(x), x \in R^m$ of degree $n$ in $m$ real varibles is called hyperbolic in the direction $e \in R^m$ 
(or $e$- hyperbolic) if for any $x \in R^m$ the polynomial 
 $p(x - \lambda e)$  in the one variable $\lambda$ has exactly $n$ real roots counting their multiplicities. We assume in this
paper that $p(e) > 0$ .
Denote an ordered vector of roots of $p(x - \lambda e)$ as 
$\lambda(x) = (\lambda_{1}(x) \geq \lambda_{2}(x) \geq ... \lambda_{n}(x)) $. It is well known that the product of roots is equal
to $p(x)$. Call $x \in R^m$ $e$-positive ($e$-nonnegative) if $\lambda_{n}(x) > 0$ ($\lambda_{n}(x) \geq 0$).
The fundamental result \cite{gar} in the theory of hyperbolic polynomials states that the set of $e$-nonnegative vectors is a 
closed convex cone. A $k$-tuple of vectors $(x_1,...x_k)$ is called $e$-positive ($e$-nonnegative) if $x_i , 1 \leq i \leq k$ are
$e$-positive ($e$-nonnegative). 
We denote the closed convex cone of $e$-nonnegative vectors as $N_{e}(p)$,
and the open convex cone of $e$-positive vectors
as $C_{e}(p)$.\\

{ \it Recent interest in the hyperbolic polynomials got sparked by the discovery \cite{gul} ,\cite{lewis} that $\log(p(x))$ is a self-concordant
barrier for the opened convex cone $C_{e}(p)$ and therefore the powerful mashinery of interior-point methods can be applied .
It is an important open problem whether this cone $C_{e}(p)$ has a semi-definite representation .}\\

It has been shown in \cite{gar} (see also \cite{khov}) that 
an $e$- hyperbolic polynomial $p$ is also
$d$- hyperbolic for all $e$-positive vectors $d \in C_{e}(p)$. \\
Let us fix $n$ real vectors $x_i \in R^m , 1 \leq i \leq n$ and define the following homogeneous polynomial:
\beqn
P_{x_1,..,x_n}(\alpha_1,...,\alpha_n) = p(\sum_{1 \leq i \leq n} \alpha_i x_i)
\eeqn
Following \cite{khov} , we define the $p$-mixed form of an $n$-vector tuple ${\bf X} = (x_1,..,x_n)$ as
\beqn
M_{p}({\bf X}) = : M_{p}(x_1,..,x_n) = \frac{\partial^n}{\partial \alpha_1...\partial \alpha_n} p(\sum_{1 \leq i \leq n} \alpha_i x_i)
\eeqn
Equivalently, the $p$-mixed form $M_{p}(x_1,..,x_n)$ can be defined by the polarization (see \cite{khov}) :
\beqn
M_{p}(x_1,..,x_n) = 2^{-n} \sum_{b_{i} \in \{-1, +1 \}, 1 \leq i \leq n }  p(\sum_{1 \leq i \leq n } b_i x_i) \prod_{1 \leq i \leq n } b_i
\eeqn

Associate with any vector $r = (r_1,...,r_n) \in I_{n,n}$ an $n$-tuple of $m$-dimensional vectors  ${\bf X}_{r}$ consisting
of  $r_i$ copies of $x_{i}  (1 \leq i \leq n) $.  
It follows from the Taylor's formula
 that
\beqn
P_{x_1,..,x_n}(\alpha_1,...,\alpha_n) = \sum_{r \in I_{n,n}} \prod_{1 \leq i \leq n } \alpha_{i}^{r_{i}} M_{p}({\bf X}_{r})
\frac{1}{\prod_{1 \leq i \leq n } r_{i}!}
\eeqn
For an $e$-nonnegative tuple ${\bf X} = (x_1,..,x_n)$, define its capacity as:
\beqn
Cap({\bf X}) = \inf_{\alpha_i > 0, \prod_{1 \leq i \leq n }\alpha _i = 1} P_{x_1,..,x_n}(\alpha_1,...,\alpha_n)
\eeqn

Probably the best known example of a hyperbolic polynomial comes from the hyperbolic geometry :
\beqn
P(\alpha_0 ,...,\alpha_k) =\alpha_0^{2} - \sum_{1 \leq i \leq k } \alpha _i^{2}  
\eeqn

This polynomial is hyperbolic in the direction $(1,0,0,...,0)$.  Another "popular" hyperbolic polynomial is
$\det(X)$ restricted on a linear real space of hermitian $n \times n$ matrices .
In this case mixed forms are just mixed discriminants , hyperbolic direction is the identity matrix $I$ ,
the corresponding closed convex cone of $I$-nonnegative vectors coincides with a closed convex cone of positive semidefinite matrices .\\ 
Less known , but very interesting , hyperbolic polynomial is the Moore determinant $Det_{(M)}(Y)$ 
restricted on a linear real space of hermitian quaternionic $n \times n$ matrices .
The Moore determinant is , essentially , the Pfaffian (see the corresponding definitions and
the theory in a very readable paper \cite{quat} ) .

We use in this paper the following class of hyperbolic in the direction $(1,1,...,1)$  polynomials of degree $k$ :\\ 
$Q(\alpha_1,...,\alpha_k) = M_{p}(\sum_{1 \leq i \leq k} \alpha_i x_i,...,\sum_{1 \leq i \leq k} \alpha_i x_i , x_{k+1},...,x_{n}) $,
where $p$ is a  $e$-hyperbolic polynomial of degree $n > k$ , $(x_1,..,x_n)$ is $e$-nonnegative tuple , and the $p$-mixed form \\
$M_{p}(\sum_{1 \leq i \leq k} x_i,...,\sum_{1 \leq i \leq k} x_i , x_{k+1},...,x_{n}) > 0$. \\ 

We make a substantial use of the following very recent result \cite{lax} , which is a rather direct corollary of
\cite{hel} , \cite{vin} . 
\thm
Consider  a homogeneous polynomial $p(y_1,y_2,y_3))$ of degree $n$ in $3$ real variables which is
hyperbolic in the direction $(0,0,1)$. Assume that $p(0,0,1)= 1$ . 
Then there exists two $n \times n$ real symmetric matrices $A,B$ such that 
$$
p(y_1,y_2,y_3)) = \det(y_1 A + y_2 B + y_3 I).
$$
\ethm 
It has been shown in \cite{hyp} that most of known facts , and some opened problems as well , 
about hyperbolic polynomials
follow from Theorem 1.5 .
\section{A hyperbolic analogue of the Rado theorem}

\dfn
Consider  a homogeneous polynomial $p(x) , x \in R^m$ of degree $n$ in $m$ real variables which is
hyperbolic in the direction $e$.Denote an ordered vector of roots of $p(x - \lambda e)$ as 
$\lambda(x) = (\lambda_{1}(x) \geq \lambda_{2}(x) \geq ... \lambda_{n}(x)) $ . We
define the $p$-rank of $x \in R^m$ in direction $e$ as $Rank_{p}(x) = |\{i : \lambda_{i}(x) \neq 0 \} |$.
It follows from Theorem 1.5 that the $p$-rank of $x \in R^m$ in any direction $d \in C_{e}$ is equal to the $p$-rank of $x \in R^m$
in direction $e$ , which we call the $p$-rank of $x \in R^m$ .
\edfn
Consider the following polynomial in one variable $D(t) = p(td+x) = \sum_{0 \leq i \leq n} c_{i} t ^{i} $.
It follows from the identity (11) that \\
\begin{eqnarray}
&c_{n} = M_{p}(d,..,d) (n!)^{-1} = p(d), \\ \nonumber
&c_{n-1} = M_{p}(x,d,..,d) (1! (n-1)!)^{-1},...,\\ \nonumber
&c_{0} = M_{p}(x,..,x) (n!)^{-1} = p(x).
\end{eqnarray}
Let $(\lambda_{1}^{(d)}(x) \geq \lambda_{2}^{(d)}(x) \geq... \geq \lambda_{n}^{(d)}(x))$ be the (real) roots of $x$
in the $e$-positive direction $d$, i.e. the roots of the equation $p(td-x) = 0$ .
Define (canonical symmetric functions) :
$$
S_{k,d}(x) = \sum_{1 \leq i_1 < i_2 <  ...< i_k \leq n} \lambda_{i_{1}}(x) \lambda_{i_{2}}(x)... \lambda_{i_{k}}(x).
$$
Then $S_{k,d}(x) = \frac{c_{n-k}}{c_{n}}$ . Clearly if $x$ is $e$-nonnegative then for any e-positive vector $d$
the $p$-rank $Rank_{p}(x) =\max \{k : S_{k,d}(x) > 0 \}$ .
The next theorem , which we prove in Appendix {\bf A} , is the main mathematical result of this paper . Our main tool is
Theorem 1.5 , which facilites a rather easy induction . We also use a particularly easy case of the Rado theorem (see
Remark A.5 for the details ) .

\thm
Consider  a homogeneous polynomial $p(x) , x \in R^m$ of degree $n$ in $m$ real variables which is
hyperbolic in the direction $e$ , $p(e) > 0$. Let  ${\bf X}=(x_1,...x_n) , x_i \in R^m$ be $e$-nonnegative $n$-tuple of $m$-dimensional vectors ,
i.e. $x_i , 1 \leq i \leq n$ are
$e$-nonnegative .\\
Then the $p$-mixed form $M_{p}({\bf X}) = : M_{p}(x_1,..,x_n)$ is positive iff the following generalized Rado conditions hold :
\beqn
Rank_{p}(\sum_{i \in S} x_i) \geq |S| \mbox{ \ for all \ } S \subset \{1,2,...,n \} . 
\eeqn
\ethm

\dfn
Call a homogeneous polynomial $p(\alpha) , \alpha \in R^n$ of degree $n$ in $n$ real variables $P$-hyperbolic
if it is hyperbolic in direction $e =(1,1,...1)$ (vector of all ones) , $p(e) > 0$ and all the canonical orts $e_i , 1 \leq i \leq n$ (rows of the identity matrix $I$ )
are $e$-nonnegative . In other words ,a homogeneous polynomial $p(\alpha) , \alpha \in R^n$ of degree $n$ in $n$ real variables is $P$-hyperbolic
if it is $e$-hyperbolic and its closed cone of $e$-nonnegative vectors contains the nonnegative orthant $R^{n}_{+} = 
\{(x_1,...,x_n) : x_i \geq 0 , 1 \leq i \leq n \}$ . It follows from \cite{khov} that the coefficients of$P$-hyperbolic 
polynomials are nonnegative real numbers .\\
(Notice that the class of $P$-hyperbolic polynomials coincides with the class of polynomials 
$P_{x_1,..,x_n}(\alpha_1,...,\alpha_n) = p(\sum_{1 \leq i \leq n} \alpha_i x_i)$ , where $p$ is $e$ -hyperbolic polynomial of
degree $n$ in $m$ real variables , a $n$-tuple $(x_1,..,x_n)$ of $m$-dimensional real vectors is $e$-nonnegative and
$\sum_{1 \leq i \leq n} x_i$ is $e$-positive . ) \\
Call a homogeneous polynomial $q(\alpha) , \alpha \in R^n$ of degree $n$ in $n$ real variables with nonnegative
coefficients $S$-hyperbolic if there exists a $P$-hyperbolic polynomial $p$ such that $supp(p) = supp(q)$ .
\edfn
(One natural class of $S$-hyperbolic polynomials ,not all of them  $P$-hyperbolic , is $Vol(\alpha_1 X_1 +...\alpha_n X_n)$ ,
where $X_i$ are convex compact subsets of $R^{n}$ . See the corresponding not $P$-hyperbolic example in \cite{khov} .)

\cor
Let  $q(\alpha) , \alpha \in R^n$ be $S$-hyperbolic polynomial of degree $n$ . \\
Then $CO(supp(q)) \cap I_{n,n} = supp(q)$ .
\ecor
\prf
It is enough to prove the corollary for $P$-hyperbolic polynomials. I.e. suppose that
$q(\alpha_1,...,\alpha_n) = p(\sum_{1 \leq i \leq n} \alpha_i x_i)$ , where $p$ is $e$ -hyperbolic polynomial of
degree $n$ in $m$ real variables , a $n$-tuple $(x_1,..,x_n)$ of $m$-dimensional real vectors is $e$-nonnegative and
$\sum_{1 \leq i \leq n} x_i$ is $e$-positive . Then $r=(r_1,r_2,...,r_n) \in supp(q)$ iff the $p$-mixed form
$M_{p}({\bf X}_{r}) > 0$ , where the $n$-tuple ${\bf X}_{r}$ consists of $r_i$ copies of $x_i ,  1 \leq i \leq n $.
Let $r^{(0)} = (r_1^{(0)},...,r_n^{(0)}) \in CO(supp(q)) .$ I.e. there exist $r^{(j)} \in supp(q) , 1 \leq j \leq n$
such that $r^{(0)} = \sum_{1 \leq j \leq n} a_j r^{(j)}$  and $a_j \geq 0 ,  \sum_{1 \leq j \leq n} a_j = 1 $ . \\
Let $r^{(j)} = (r_{1}^{(j)} , ... ,r_{n}^{(j)} ) ,  0 \leq j \leq n$ . As $r^{(j)} \in supp(q) , 1 \leq j \leq n$ 
thus $M_{p}({\bf X}_{r^{(j)}}) > 0 ,1 \leq j \leq n$ . It follows from Theorem 2.2 (only if part ) that
$$
Rank_{p}(\sum_{i \in S} x_i) \geq  \sum_{i \in S}  r_{i}^{(j)} \mbox{ \ for all \ } S \subset \{1,2,...,n \} ; 1 \leq j \leq n .
$$
Therefore
$$
Rank_{p}(\sum_{i \in S} x_i) \geq \sum_{i \in S} \sum_{1 \leq j \leq n} a_j r_{i}^{(j)} = \sum_{i \in S}  r_{i}^{(j)} , S \subset \{1,2,...,n \} .
$$
Using the "if" part of Theorem 2.2 we get that  $M_{p}({\bf X}_{r^{(0)}}) > 0$ and thus $r^{(0)} \in supp(q)$ .
\eprf

\cor
Let  $q(x) , x \in R^n$ be $S$-hyperbolic polynomial of degree $n$ . Then the following conditions are equivalent
\begin{enumerate}
\item
$e \in CO(supp(q))$ .
\item
$e \in supp(q)$ , i.e. $\frac{\partial^n}{\partial \alpha_1...\partial \alpha_n} q(x) > 0 $ .
\item
$Cap(p) = : \inf_{\alpha_i > 0, \prod_{1 \leq i \leq n }\alpha _i = 1} q(\alpha_1,...,\alpha_n) > 0 $.
\item
For all $\epsilon > 0$ there exists a vector $(\alpha_1,...,\alpha_n)$ with positive entries
such that the following inequality holds :
\beqn
\sum_{1 \leq i \leq n} |\frac{ \alpha_{i } \frac{\partial}{\partial \alpha_{i } } q(\alpha_1,...,\alpha_n) } {q(\alpha_1,...,\alpha_n)} - 1|^{2}
\leq \epsilon .
\eeqn
\item
There exists a vector $(\alpha_1,...,\alpha_n)$ with positive entries
such that the following inequality holds :
\beqn
\sum_{1 \leq i \leq n} |\frac{ \alpha_{i } \frac{\partial}{\partial \alpha_{i } } q(\alpha_1,...,\alpha_n) } {q(\alpha_1,...,\alpha_n)} - 1|^{2}
\leq \frac{1}{n} .
\eeqn
\item
For all subsets $S \subset \{1,2,...,n \}$ the following inequality holds :
\beqn
\sum _{i \in S} r_i  \geq |S| \mbox{\ for all \ } (r_1,...,r_n) \in supp(q) .
\eeqn 

\end{enumerate}
\ecor
(We sketch a proof in Appendix {\bf C} . ) \\

The following result , which we prove in Appendix {\bf B} , is a "polynomial" generalization of Lemma 4.2 in \cite{lsw} .
\pro
The condition (17) implies the condition (18) for all homogeneous polynomial $q(x) , x \in R^n$ of degree $n$ in $n$ real variables with nonnegative
coefficients . 
\epro
\subsection{Doubly-stochastic polynomials}
Inequalities (16), (17) above suggest the following definition .
\dfn
A homogeneous polynomial $q(x_1,...,x_n)$ of degree $n$ in $n$  variables is called doubly-stochastic
if its coefficients are nonnegative real numbers and $\frac{\partial}{\partial x_{i } } q(1,1,...,1) = 1$ for all $1 \leq i \leq n$.
The doubly-stochastic defect of the polynomial $q$ is defined as $DS(q) = \sum_{1 \leq i \leq n} (\frac{\partial}{\partial x_{i } } q(1,1,...,1) - 1)^{2}$
\edfn
\lem
\begin{enumerate}
\item A homogeneous polynomial $q(x_1,...,x_n)$ of degree $n$ in $n$  variables with nonnegative real coefficients is doubly-stochastic iff
$q(1,1,...,1)=1$ and $q(x_1,...,x_n) \geq \prod_{1 \leq i \leq n} x_i$ for all real vectors $(x_1,...,x_n) \in R^n$ with positive coordinates 
(in other words if $q(1,1,...,1)=1$ and $Cap(q)= 1$ ).
\item A homogeneous polynomial $q(x_1,...,x_n)$ of degree $n$ in $n$  variables is $P$-hyperbolic and doubly-stochastic iff
$q(1,1,...,1)=1$ and $|q(z_1,...,z_n)| \geq \prod_{1 \leq i \leq n} Re(z_i)$ for all complex vectors vectors $(z_1,...,z_n) \in C^n$ with positive 
real parts .
\item
If a sequence $q_{i}$ of homogeneous polynomials of degree $n$ in $n$  variables with nonnegative real coefficients converges to
a doubly-stochastic polynomial then $\lim_{i \rightarrow \infty} Cap(q_{i}) = 1$.
\item
The capacity $Cap(q)$ is a continuous functional (but not even Lipshitz) on a convex closed cone of homogeneous polynomials of degree $n$ in $n$  variables with nonnegative real coefficients .
\item
If $q$ is homogeneous polynomial of degree $n$ in $n$  variables with nonnegative real coefficients then $Cap(q) > 0$ iff
the exists a sequence $X_j = (x_{1,j},...,x_{n,j})$ of vectors with positive real coordinates such that
$$
\sum_{1 \leq i \leq n} |\frac{ x_{i} \frac{\partial}{\partial x_{i} } q(x_{1,j},...,x_{n,j}) } {q(x_{1,j},...,x_{n,j})} - 1|^{2} \rightarrow 0
$$
And in this case $Cap(q) = \lim_{j \rightarrow \infty} \frac{q(X_{j})}{\prod_{1 \leq i \leq n} x_{i,j}}$.
\end{enumerate} 
\elem
\xam
A multilinear polynomial \\ 
$q(x_1,...,x_n) = \prod_{1 \leq i \leq n} \sum_{1 \leq j \leq n} a(i,j) x_j$ is doubly-stochastic and $P$-hyperbolic iff \\
the square matrix $B =\{\frac{a(i,j)}{\sum_{1 \leq k \leq n} a(i,k)} : 1 \leq i,j \leq n \}$ is doubly stochastic in the standard
meaning .
\exam

The next theorem is another Corollary of Theorem 2.2 ; its main point is in introducing a hyperbolic analog of Van der Waerden conjecture .
\thm
\begin{enumerate}
\item Consider the set $PHDS(n)$ of all $P$-hyperbolic doubly-stochastic homogeneous polynomials $q(x_1,...,x_n)$ of degree $n$ in $n$ variables.
The set $PHDS(n)$ is a compact and the following inequality holds
\begin{eqnarray*}
&\inf_{q \in PHDS(n)}\frac{\partial^n}{\partial x_1...\partial x_n} q(x_1,...,x_n)  = \\
&\min_{q \in PHDS(n)}\frac{\partial^n}{\partial x_1...\partial x_n} q(x_1,...,x_n) =:VdW(n) > 0
\end{eqnarray*}
\item For any $P$-hyperbolic homogeneous polynomials $q(x_1,...,x_n)$ of degree $n$ in $n$ variables the following inequalities hold
$$
VdW(n) \leq \frac{\frac{\partial^n}{\partial x_1...\partial x_n} q(x_1,...,x_n)}{Cap(q)} \leq 1 .
$$
\end{enumerate}
\ethm
\cnj {\bf Hyperbolic Van der Waerden conjecture }\\ 
$$
VdW(n) = \frac{n!}{n^{n}} ?
$$
\ecnj

The next result , a direct corollary of Lemma 2.10 in \cite{hyp} , is a generalization of Proposition 4.2 in \cite{GS1}.

\lem
Let $q$ be $P$-hyperbolic homogeneous polynomial of degree $n$ in $n$ variables . Then the following inequalities hold
for all vectors $(x_1,...,x_n)$ with positive real coordinates.
\begin{enumerate}
\item
\beqn
Cap(q) \leq \frac{q(x_1,...,x_n)}{\prod_{1 \leq i \leq n} x_i} \prod _{1 \leq i \leq n} x_i \frac{\frac{\partial}{\partial x_{i} } q(x_1,...,x_n) } {q(x_1,...,x_n)} 
\eeqn
\item If $\log(\frac{q(x_1,...,x_n)}{\prod_{1 \leq i \leq n} x_i}) - \log(Cap(q))  \leq \frac{\epsilon}{10}$ with $0 \leq \epsilon \leq 1$ then
\beqn
\sum_{1 \leq i \leq n} |\frac{ x_{i } \frac{\partial}{\partial x_{i } } q(x_1,...,x_n) } {q(x_1,...,x_n)} - 1|^{2}
\leq \epsilon
\eeqn
\end{enumerate}
\elem

\xam
Consider the following homogeneous polynomial of degree $n$ in $n$ variables : $p(x_1,...,x_n) = \sum_{1 \leq i \leq n} x_i^{n}$.
Then $Cap(p) = n$ and 
\begin{eqnarray*}
&\frac{p(x_1,...,x_n)}{\prod_{1 \leq i \leq n} x_i} \prod _{1 \leq i \leq n} x_i \frac{\frac{\partial}{\partial x_{i} } p(x_1,...,x_n) } {p(x_1,...,x_n)} = \\
&= n^{n} (\frac{ \prod _{1 \leq i \leq n} x_i} {\sum_{1 \leq i \leq n} x_i^{n}})^{n-1} \leq n = Cap(p).
\end{eqnarray*}
The moral of this example is that the inequality (19) does not hold for all homogeneous polynomials with nonnegative coefficients ,
this inequality is a nontrivial necessary condition for the \\
$P$-hyperbolicity . It is interesting to notice that the inequality (19) implies the determinantal Hadamard inequality .
\exam
\rem
Let 
$$
p(x_1,...,x_n) =\sum_{ (r_1,...,r_n) \in I_{n,n} } a_{(r_1,...,r_n) } \prod_{1 \leq i \leq n} x_{i}^{r_{i}}
$$ 
be homogeneous polynomial of degree $n$ in $n$ real variables. Assume that its coefficients are nonnegative and sum to one , i.e. that $p(1,1,...,1)=1$.
Associate with this polynomial a random integer vector 
$$
Z_{p} = (z_1,...,z_n) \in  I_{n,n} : Prob\{Z_{p} =(r_1,...,r_n)\}= a_{(r_1,...,r_n) }. 
$$
Then 
$$
E(Z_{p}) = (\frac{\partial}{\partial x_{1} } p(1,1,...,1),..., \frac{\partial}{\partial x_{n} }p(1,1,...,1)).
$$
Therefore $p$ is doubly-stochastic iff $E(Z_{p}) =(1,1,...,1)$ ; there exists a doubly-stochastic polynomial $q$ such that $supp(q) \in supp(p)$
iff $(1,1,...,1) \in CO(supp(p))$. \\
It follows from Corollary 2.5 that if $p$ is doubly-stochastic and $S$-hyperbolic then
$Prob\{Z_{p} = E(Z_{p}) \} > 0$ . And the hyperbolic van der Waerden conjecture can be reformulated as : \\
{\it If $p$ is doubly-stochastic and $P$-hyperbolic then \\
$Prob\{||Z_{p} = E(Z_{p})|| < \sqrt{2} \} \geq \frac{n!}{n^{n}}$ } . \\
Perhaps some kind of the {\ measure concetration } is present here ?
\erem

\rem
The problem to find out a positive real solution of the inequality (20) is a far reaching
generalization of scaling of matrices with nonnegative entries (the corresponding
polynomials are multilinear )\cite{lsw} ,\cite{nr} and scaling
of tuples of PSD matrices (the corresponding polynomials are determinantal )
 \cite{GS} , \cite{GS1} . Part 2 of Lemma 2.12 allows
to generalize results of \cite{nr} , \cite{GS} , \cite{GS1} to $P$-hyperbolic polynomials ,
even in the black-box setting . Can it be done for all homogeneous
polynomials with , say , integer nonnegative coefficients ?
\erem

\section{The ellipsoid algorithm }
Consider a homogeneous polynomial $q(x) , x \in R^n$ of degree $n$ in $n$ real variables with nonnegative integer
coefficients . Associate with such $q$ the following convex functional 
$$
f(y_1,...,y_n) = \log(q(e^{y_{1}},e^{y_{2}},...,e^{y_{n}})) .
$$

\pro
The following conditions are equivalent
\begin{enumerate}
\item
$e=(1,1,..,1) \in CO(supp(q))$ .
\item
$ \inf_{y_1+...+y_n =0} f(y_1,...,y_n) \geq 0 $.
\item
If $e = (1,1,..,1) \notin CO(supp(q))$ then \\
$ \inf_{y_1+...+y_n =0} f(y_1,...,y_n) = - \infty $. \\
Let $dist(e, CO(supp(q))) = \Delta^{-1} > 0$ and $Q = \log(q(e))$ . Define $\gamma = (Q+1)\Delta $ . 
Then
\begin{eqnarray*}
&\inf_{y_1+...+y_n =0 , (|y_1|^{2}+...+|y_n|^{2})^{\frac{1}{2}} \leq \gamma  } f(y_1,...,y_n) = \\
&\min_{y_1+...+y_n =0 , |y_1|^{2}+...+|y_n|^{2} \leq \gamma} f(y_1,...,y_n) \leq -1 .
\end{eqnarray*}
\end{enumerate}
\epro

\prf
Our proof is a straigthforward application of the concavity of the logarithm on the positive semi-axis and
of the Hanh-Banach separation theorem . It will be included in the full version .
\eprf

Proposition 3.1 suggests the following natural approach to solve Problem 2 , i.e. to decide
whether $e=(1,1,..,1) \in CO(supp(q))$ or not : \\
find $\min_{y_1+...+y_n =0 , |y_1|^{2}+...+|y_n|^{2} \leq \gamma } f(y_1,...,y_n)$ with absolute accuracy
$ \frac{1}{3}$ . If the resulting value is greater than or equal $-\frac{1}{3}$ then $e=(1,1,..,1) \in CO(supp(q))$ ;
if the resulting value is less than or equal $-\frac{2}{3}$  then $e=(1,1,..,1) \notin CO(supp(q))$ .
And , of course , it is natural to use the ellipsoid method .
Our main tool is the following property of the ellipsoid algorithm
\cite{ne:ne}: 
For a prescribed accuracy $\delta > 0$, 
it finds a $\delta$-minimizer of a differentiable convex function $f$ in a ball
$B$, that is a point $x_{\delta} \in B$ with $f(x_{\delta}) \le
\min_B f + \delta$, in no more than
$$
O\left(n^2 \ln \left(\frac{2\delta +
\mbox{Var}_B(f)}{\delta}\right)\right),~~~~~~~~(\mbox{Var}_B(f) =
\max_B f - \min_B f )
\label{ell_cost}
$$
iterations. Each iteration requires a single
computation of the value and of the gradient of $f$ at a given point,
plus $O(n^2)$ elementary operations to run the algorithm itself. In
our case, this is easily seen to cost at most $O(n^2)$ oracle calls and
$O(n)$ elementary arithmetic operations . \\
We have the $n-1$-dimensional ball $B_{\gamma} = \{(y_1,...,y_n) : y_1+...+y_n =0 , |y_1|^{2}+...+|y_n|^{2} \leq \gamma \}$.
A straigthforward computations show that 
$$
Var_B(f) \leq \log(q(1,1,..,1) e^{\gamma n} ) - \log(q(1,1,..,1) e^{-\gamma n} ) \leq 2 \gamma n ,
$$
giving that $O(n^{2} ( \ln(n) + \ln(\gamma ) )$ iterations of the ellipsoid method needed to solve Problem 2 ,
it amounts to $O(n^{4} ( \ln(n) + \ln(\gamma ) )$ oracle calls .  The quantity $O(n^{4} ( \ln(n) + \ln(\gamma ) )$ is
polynomial in $n$ even if $\gamma$ is exponentially large ($dist(e, CO(supp(q)))$ is  exponentially small ).
The problem is that if  $\gamma$ is exponentially large ( which can happen ) then we need to call oracles on inputs
with exponential bit-size . \\
Putting things together , we get the following conclusion : \\
{\it If it is promised that either $e=(1,1,..,1) \in CO(supp(q))$ or $dist(e, CO(supp(q))) \geq poly(n)^{-1}$ for
some fixed polynomial $poly(n)$ then Problem 2 can be solved by a deterministic polynomial-time oracle algorithm
based on the ellipsoid method .} \\
And at this point we can say nothing about Problem 1 , i.e. deciding whether $e=(1,1,..,1) \in supp(q)$ or not .
Corollary 2.5 says that if $q$ is $S$-hyperbolic polynomial then Problem 1 and Problem 2 are equivalent ;
moreover if $e=(1,1,..,1) \notin supp(q)$ then here exists nonempty $S \subset \{1,2,...,n \}$ such that
\beqn
\sum _{1 \leq i \leq n} r_i s_i < \sum _{1 \leq i \leq n}  s_i = |S| \mbox{\ for all} (r_1,...,r_n) \in supp(q) ,
\eeqn ,
where $(s_1,...,s_n)$ is a characteristic function of the subset $S$ , i.e. $s_i = 1$ if $i \in S$ , and $s_i = 0$ otherwise . \\
Notice that if (21) holds then the distance $dist(e,CO(supp(q)))$  from the vector $e = (1,...,1)$ to the Newton polytope $CO(supp(q))$
is at least $\sqrt{\frac{n}{|S| (n- |S|) }} \geq \frac{2}{\sqrt{n}}$ . Thus we have the next theorem .
\thm
Problem 1 and Problem 2 are equivalent for $S$-hyperbolic polynomials .\\
There exists a deterministic polynomial-time oracle algorithm
solving Problem 1 for a given $S$-hyperbolic polynomial $q(\alpha_1,...,\alpha_n)$ with integer coefficients .\\ 
It requires 
$O(n^{4} (\ln(n) + \ln(\ln(q(1,1,...,1)))$ oracle calls and it bit-wise complexity (which is roughly the radius
of the ball $B_{\gamma}$ )  is $O(n^{\frac{1}{2}} \ln(q(1,1,...,1)))$ .
\ethm

\section{Hyperbolic Sinkhorn scaling}
We will discuss briefly in this section another method , which is essentially a large step version of gradient descent .
\dfn
Consider an $e$-nonnegative tuple ${\bf X} = (x_1,..,x_n)$ such that the sum of its components
 $S({\bf X}) = d = \sum_{1 \leq i \leq k} x_i$ is $e$-positive. Define $tr_{d}(x)$ as a sum of roots
of the univariate polynomial equation $p(x -t d)= 0$.\\
Define the following map (Hyperbolic Sinkhorn Scaling) acting on such tuples:
$$
HS({\bf X}) = {\bf Y} = (\frac{x_{1}}{tr_{d}(x_{1})},..., \frac{x_{n}}{tr_{d}(x_{n})})
$$
Hyperbolic Sinkhorn Iteration ({\bf HSI}) is the following recursive procedure:
\begin{eqnarray*}
&{\bf X}_{j+1} = HS({\bf X}_{j}), j \geq 0, ~{\bf X}_{0} \mbox{ is an $e$-nonnegative tuple with }\\
& \sum_{1 \leq i \leq k} x_i \in C_{e}\;.
\end{eqnarray*}
We also define the doubly-stochastic defect of $e$-nonnegative tuples with $e$-positive sums as
$$
DS({\bf X}) = \sum_{1 \leq i \leq k} (tr_{d}(x_{i}) - 1)^{2} ; \sum_{1 \leq i \leq k} x_i = d \in C_{e}
$$
\edfn

We can define the map $HS(.)$ directly in terms of the $P$-hyperbolic polynomial 
$$
Q(\alpha_1,...,\alpha_n) = P_{x_1,..,x_n}(\alpha_1,...,\alpha_n) = p(\sum_{1 \leq i \leq n} \alpha_i x_i).
$$
Indeed, if $ \sum_{1 \leq i \leq n} \alpha_i x_i =d \in C_{e}$ then
\beqn
tr_{d}(\alpha_i x_i) = \frac{ \alpha_{i } \frac{\partial}{\partial \alpha_{i } } Q(\alpha_1,...,\alpha_n) } {Q(\alpha_1,...,\alpha_n)}
\eeqn
This gives the following way to redefine the map $HS({\bf X})$ :
$$
HS(\alpha_1,...,\alpha_n) = (\frac{Q(\alpha_1,...,\alpha_n)}{ \frac{\partial}{\partial \alpha_{1} } Q(\alpha_1,...,\alpha_n) } ,...,
\frac{Q(\alpha_1,...,\alpha_n)}{ \frac{\partial}{\partial \alpha_{n } } Q(\alpha_1,...,\alpha_n) } ) ,
$$
for $\alpha_{i } > 0 , 1 \leq i \leq n $. \\
And correspondingly the doubly-stochastic defect of $(\alpha_1,...,\alpha_n)$ is equal to
$$
\sum_{1 \leq i \leq n} |\frac{ \alpha_{i } \frac{\partial}{\partial \alpha_{i } } Q(\alpha_1,...,\alpha_n) } {Q(\alpha_1,...,\alpha_n)} - 1|^{2} ,
$$
the same as the left side of (17 ) . Notice that $\sum_{1 \leq i \leq n} tr_{d}(x_{i}) = n $ by the  Euler's identity .
\xam
Consider the following hyperbolic polynomial in $n$ variables: $p(z_1,...,z_n) = \prod_{1 \leq i \leq n} z_i$.
It is $e$- hyperbolic for $e = (1,1,...,1)$. And $N_{e}$ is a nonnegative orthant, $C_{e}$ is a positive orthant.
An $e$-nonnegative tuple ${\bf X} = (x_1,..,x_n)$ can be represented by an $n \times n$ matrix $A_{{\bf X}}$ with nonnegative entries:
the $i$th column of $A$ is a vector $x_i \in R^{n}$. If $Z = (z_1,...,z_n) \in R^{n}$ and
$d = (d_1,...,d_n) \in R^{n} ; z_i > 0, 1 \leq i \leq n$, then $tr_{d}(Z) = \sum_{1 \leq i \leq n} \frac{z_{i}}{d_{i}}$.
Recall that for a square matrix $A = \{a_{ij}: 1 \leq i,j \leq N\}$ row
scaling is
defined as
$$
R(A) = \{ \frac{a_{ij}}{\sum_j a_{ij}} \}, 
$$
column scaling as $C(A) = \{ \frac{a_{ij}}{\sum_i a_{ij}} \}$
assuming that all
denominators are nonzero.
The iterative process $...CRCR(A)$ is called {\em Sinkhorn's iterative scaling} (SI).
In terms of the matrix $A_{{\bf X}}$ the map $HS({\bf X})$ can be realized as follows:
$$
A_{HS({\bf X})} = C(R(A_{{\bf X}}))
$$
So, the map $HS({\bf X})$  is indeed a (rather far-reaching) generalization of Sinkhorn's scaling. Other generalizations (not all hyperbolic)
can be found in \cite{GY}, \cite{stoc}, \cite{arxiv}.
\exam
Lemma 2.10 from \cite{hyp}  allows to use ({\bf HSI}) to solve Problem 1 for $P$-hyperbolic polynomials $q$ in the same way as it was done for the perfect matching problem in \cite{GY} , \cite{lsw} ; and for the Edmonds' problem in \cite{stoc} . 
The corresponding complexity is $O(n \log(q(e)))$ iterations of ({\bf HSI}) , which can be done in  $O(n^{3} \log(q(e)))$ oracle calls .
The algorithm works in the following way : \\
{\it Run $K =O(n \log(q(e)))$ Hyperbolic Sinkhorn Iterations  ${\bf X}_{j+1} = HS({\bf X}_{j}) $   ; if $DS({\bf X}_{i}) \leq \frac{1}{n}$ for
some $ i \leq K$ then the $p$-mixed form $M_{p}({\bf X}_{0}) > 0$ , and $M_{p}({\bf X}_{0})= 0$ otherwise .}

\section{ Half-Plane Property }
The following definition is from \cite{half}.
\dfn
A polynomial $P(z_1,...,z_n)$ in $n$ complex variables is said to have the "half-plane property"
if $P(z_1,...,z_n) \neq 0$ provided $Re(z_i) > 0$ .
\edfn
In a control theory literature (see \cite{khar} ) the same property is called {\it Wide sense stability} .
And {\it Strict sense stability} means that \\
$P(z_1,...,z_n) \neq 0$ provided $Re(z_i) \geq 0$ .

The following simple fact shows that for homogeneous polynomials the "half-plane property" is , up to a
single factor , the same as $P$-hyperbolicity .
\fac
A homogeneous polynomial $R(z_1,...,z_n)$ has the "half-plane" property if and only if
the exists real $\alpha$ such that the polynomial $e^{i\alpha} R(z_1,...,z_n)$ is $P$-hyperbolic polynomial
with real nonnegative coefficients  .
\efac
\prf
\begin{enumerate}
\item
Suppose that $R(z_1,...,z_n) = e^{-i\alpha}Q(z_1,...,z_n)$ where $\alpha$ is real and $Q$ is $P$-hyperbolic.
Then $Q$ is $(1,1,...,)$-hyperbolic and all real vectors $(x_1,...,x_n)$ with positive coordinates are $(1,1,...,)$-positive .
Therefore $Q$ is $ (x_1,...,x_n)$-hyperbolic for all  real vectors $(x_1,...,x_n) \in R^{n}_{++}$ with positive coordinates .It follows
that $|R(x_1 + i y_1,...,x_n + i y_n)| = |Q(x_1 + i y_1,...,x_n + i y_n)| = |Q(x_1,...,x_n) \prod_{1 \leq k \leq n} (1 + i \lambda_{k})|$ ,
where $(\lambda_{1},...,\lambda_{n})$ are real roots of the real vector $(y_1,...,y_n)$ in the direction $(x_1,...,x_n)$.\\
This gives the following inequality , which is  equivalent to  the "half-plane property" of $R$ :
\begin{eqnarray}
&|R(x_1 + i y_1,...,x_n + i y_n)| \geq |R(x_1,...,x_n)| = \\ \nonumber
&= |Q(x_1,...,x_n)| > 0 : \\ \nonumber
&(x_1,...,x_n) \in R^{n}_{++} ,(y_1,...,y_n) \in R^{n} 
\end{eqnarray}
\item
Suppose that $R(z_1,...,z_n)$ has the "half-plane property" and consider the roots of the following polynomial equation in one complex variable :
$P(x_1 - z ,x_2 -z,...,x_n - z) = 0$ , where $(x_1,...,x_n) \in R^{n}$ is a real vector , $z = x +iy \in C$. If the imaginery part $Im(z) =y$
is not zero then , using the homogeniuty , $R( i \frac{x - x_{1}}{y} + 1,...,i \frac{x - x_{n}}{y} + 1) = 0$ , which is impossible as
$R$ has the "half-plane property". Therefore all roots of $R(X -t e) =0$ are real for all real vectors $X \in R^{n}$ (here $e = (1,1,...,1)$).
In the same way all roots of $R(X -t e) =0$ are real positive numbers if $X \in R^{n}_{++}$ . It follows that if $X \in R^{n}$ then
$R(X) = R(e) \prod{1 \leq k \leq n} \lambda_{k}(X)$ ,
where $(\lambda_{1},...,\lambda_{n})$ are ( real ) roots of the equation $R(X -t e) =0$ . Thus the polynomial $(\frac{1}{R(e)}) R$ takes
real values on $R^n$ and therefore its coefficients are real . In other words , the polynomial $(\frac{1}{R(e)}) R$ is $P$-hyperbolic .
If $R(1,1,...,1) = e^{- i \alpha} |R(1,1,...,1)|$ then the polynomial $e^{i\alpha} R$ is also $P$-hyperbolic .\\
(Recall that the coefficients of any $P$-hyperbolic polynomial $p$ are nonnegative for they are $p$-mixed forms of $e$-nonnegative tuples ,
and $p$-mixed forms of $e$-nonnegative tuples are nonnegative if $p(e) > 0$ \cite{khov}.)
\end{enumerate}
\eprf

We use this observation to show that Theorem 2.2 in this paper implies (and seriously strengthens)
Theorem 7.2 in \cite{half} , which is the main result of a very long recent paper \cite{half} . 
\subsection{Submodularity and hyperbolicity}
Let $p$ be a $P$-hyperbolic polynomial of degree $n$ in $n$ variables .
It follows from Theorem 2.2 that $r = (r_1,r_2,...,r_n) \in supp(p)$ if and only if the following
inequalities hold  :
$$
r(S) = \sum_{i \in S} r_i \leq R(S) = Rank_{p}(\sum_{i \in S} e_i) ; S \subset \{1,2,...,n \}.
$$

\fac
The functional $R(S)= Rank_{p}(\sum_{i \in S} e_i)$ is normalized , i.e. $R(\emptyset) = 0$ ,
and submodular , i.e. $R(A \cup B) \leq R(A) + R(B) - R(A \cap B) : A,B \subset \{1,2,...,n \}$ .
\efac
\prf
Associate with two subsets $A,B \subset \{1,2,...,n \}$ the following three $e$-nonnegative vectors : \\
$ x = \sum_{i \in A \setminus (A \cap B)} e_i , y = \sum_{i \in A \cap B } e_i ,  z = \sum_{i \in B \setminus (A \cap B)} e_i . $ \\
We need to prove the inequality $Rank_{p}(x+y+z) \leq Rank_{p}(x) + Rank_{p}(z) - Rank_{p}(y)$.
This inequality is obvious and well known for positive semidefinite matrices .
The extension to $e$-nonnegative vectors respect to $e$-hyperbolic polynomial $p$ is
done in the same way as in the proof of Corollary A.3 :
consider a hyperbolic in the direction $(1,1,1)$ polynomial 
$$
L(\alpha_1,\alpha_2,\alpha_3) = M_{p}(k,...,k,e,...,e) , k = \alpha_1 x + \alpha_2 y + \alpha_3 z ;
$$
where the vectors $x,y,z$ are $e$-nonnegative respect to hyperbolic polynomial $p$ ,
and the tuple $(k,...,k,e,...,e)$ consists of $Rank_{p}(x+y+z)$ copies of $k$ and
$n- Rank_{p}(x+y+z)$ copies of $e$. After that apply Theorem 1.5 .
\eprf

\cor
\begin{enumerate}
\item
A support $supp(p)$ of $P$-hyperbolic polynomial $p$ is an intersection of the integral polymatroid
$\{(r_1,...,r_n) : r(S) = \sum_{i \in S} r_i \leq R(S) = Rank_{p}(\sum_{i \in S} e_i) ; S \subset \{1,2,...,n \} \}$
with the hyperplane $\{(r_1,...,r_n) : \sum_{1 \leq i \leq n} r_i = n \}$.
\item
A support $supp(R)$ of any  polynomial $R$ with the "half-plane" property is a {\it jump system} .
\end{enumerate}
\ecor
\prf
(Consult  \cite{jump} for a definition and some properties of jump systems and integral polymatroids ) .
This Corollary follows directly Theorem 2.2 , Fact 5.3 and Proposition (3.1) in \cite{jump} .
\eprf
\ \\ 
It is quite amazing how the two communities , "hyperbolic" and "half-plane" , were not
aware about each other results for a long , long time . (Interestingly ,
two authors of \cite{half} and one author of \cite{lax} were with the same department
until very recently . Perhaps , one needs to be a dilettante to notice a bridge .)

\section{Conclusion and Acknowledgments}
Univariate polynomials with real roots appear quite often in modern combinatorics , especially in the context
of integer polytopes . We discovered in this
paper rather unexpected and very likely far-reaching connections between hyperbolic polynomials and many classical combinatorial and
algorithmic problems . (The author taught about "On hyperbolic nature of perfect marriages" as a title of this paper ,
but with the current climate it could be understood in many ways .) 
There are still several open problems . The most interesting is  Conjecture 2.11 in this paper , which is a  generalization
of the van der Waerden conjecture for permanents of doubly stochastic matrices and many others related questions .\\
For a hyperbolic in direction $(1,1,..,1)$ polynomial $Mul(y_1,...,y_n) = y_1 y_2 ...y_n$ Conjecture 2.11 is equivalent to 
the famous van der Waerden conjecture for permanents of doubly stochastic matrices , proved in \cite{ego} , \cite{fal} .
For a hyperbolic in direction $I$ polynomial $\det(X)$ , $X$ is $n \times n$ hermitian matrix , it is equivalent
to Bapat's conjecture \cite{bapat} (it was also hinted in \cite{ego} ) , proved by the author in \cite{GS} , \cite{gur} . 
It also holds for the Moore determinant $Det{(M)}(Y)$ , $Y$ is $n \times n$ quaternionic hermitian matrix ,
with the proof essentially the same as in \cite{gur} . \\
Another , equivalent form of "hyperbolic" (or "half-plane" ) van der Waerden conjecture can be formulated as follows :
\cnj
Consider a homogeneous polynomial $p(z_1,...,z_n)$ of degree $n$ in $n$ complex variables .
Assume that this polynomial satisfies the  property : \\

$|p(z_1,...,z_n)| \geq \prod_{1 \leq i \leq n} Re(z_i)$ on the domain $\{(z_1,...,z_n) :  Re(z_i) \geq 0 , 1 \leq i \leq n \}$ . \\

Is it true that $|\frac{\partial^n}{\partial z_1...\partial z_n} p | \geq \frac{n!}{n^n}$ ? . \\
(Notice that Theorem 2.10 and Fact 5.2 imply that $\frac{\partial^n}{\partial z_1...\partial z_n} p \neq 0$ .)
\ecnj
It would be very interesting and enlighting to prove Conjecture 2.11 using methods
of the theory of functions of many complex variables. Fact 5.2 , together with other results in this paper ,
makes a connection between the Complexity Theory and the theory of linear time-miltidimensional systems :
all "hard" instances of Problem 1 are necessary unstable polynomials. \\ 
Another interesting conjecture is related to the majorization :
\cnj
Consider the doubly-stochastic and $P$-hyperbolic homogeneous polynomial $p(x_1,...,x_n)$ of degree $n$ in $n$ real variables .\\
Let $\Lambda(X) \in R^n$ be a real $n$-dimensional vector , whose coordinates are the roots of the equation $p(X-te)=0$ ,
where $X \in R^n$ and $e$ is the vector of all ones . Then there exists a $n \times n$ doubly stochastic matrix $A$ such
that $\Lambda(X) = A X$.\\

(Some partial and related results in this direction can be found in \cite{hyp} ; this conjecture is true for
determinantal polynomials .)
\ecnj
A natural extension of {\bf Problem 1} is for $P$-hyperbolic polynomials
to approximate $\frac{\partial^n}{\partial x_1...\partial x_n} p(x_1,...,x_n)$
within a multiplicative factor using deterministic ( or randomized ) polynomial-time oracle algorithms .
It is not clear to the author whether known recent randomized algorithms for $(1+ \epsilon)$ approximation
of the permanent $Per(B)$ of entry-wise nonnegative matrix $B$ can be done in the "oracle fashion" ,
i.e. using only some outputs of the multilinear polynomial 
$$
q(x_1,...,x_n) = \prod_{1 \leq i \leq n} \sum_{1 \leq j \leq n} B(i,j) x_j .
$$
If hyperbolic van der Waerden conjecture is true then the technique in this paper ,
similarly to \cite{lsw} and \cite{GS} , \cite{GS1} , would produce 
a deterministic polynomial-time oracle algorithm with $\frac{n^{n}}{n!}$ multiplicative factor .\\

The technique developed in this paper can be applied to other "noble" desicion problems . For instance ,
checking factorizability of $P$-hyperbolic polynomials can be also done in deterministic oracle polynomial time .
The factorizability is closely related to the hyperbolic generalization of the indecomposability of matrix tuples \cite{GS1}. \\
This paper is probably the first one which uses Theorem 1.5 in the combinatorial context . We expect many more
such applications of Theorem 1.5 . This (very nontrivial) theorem , when in good hands , is a powerful tool allowing
reasonably simple and short proofs .\\

I would like to acknowledge a great influence of amazingly clear paper \cite{khov} .
It is my pleasure to thank Adrian Lewis for numerous as e-mail as well phone communications.
Many thanks to the fantastic library of Los Alamos National 
Laboratory and Google.

\appendix
\section{ Proof of the (main ) Theorem 2.2 }
Before proving Theorem 2.2 , we will recall some basic properties of $p$-mixed forms and prove a few auxillary results .
The following fact was proved in \cite{khov}
\fac
Consider  a homogeneous polynomial $p(x) , x \in R^m$ of degree $n$ in $m$ real variables which is
hyperbolic in the direction $e$. Then the following properties hold .
\begin{enumerate}
\item
The  $p$-mixed form $M_{p}(x_1,..,x_n)$ is linear in each $x_i , 1 \leq i \leq n $.
\item
If $x_1,x_2,..,x_{n-1}$ are $e$-nonnegative then the linear functional $l(x) = M_{p}(x_1,..,x_{n-1},x)$ is
nonnegative on the closed cone $N_{e}$ of $e$-nonnegative vectors .
\item
If the tuples $(x_1,..,x_n) , (y_1,..,y_n) , (x_1-y_1,..,x_n-y_n)$ are $e$-nonnegative then
$$
0 \leq M_{p}(y_1,..,y_n) \leq M_{p}(x_1,..,x_n).
$$
\item
Fix $e$-positive vector $d$ and consider the following  homogeneous polynomial $p_{d}(x) , x \in R^m$ of degree $n-1$ in $m$ real variables :
$p_{d}(x) =:  M_{p}(x,x,...,x,d)$ . Then $ p_{d}$ is hyperbolic in any $e$-positive direction $v \in C_{e}(p)$ .
If $g \in C_{e}(p)$ ( $e$-positive respect to the polynomial $p$ ) then also $q \in C_{v}(p_{d})$ for all
$v \in C_{e}(p)$ .
\end{enumerate}
\efac
The next fact is well known .
\fac
Consider a sequence of univariate polynomials of the same degree $n$ : $P_{k}(t) = \sum_{0 \leq i \leq n} a_{i,k} t^{i}$ .
suppose that $ \lim_{k \rightarrow \infty} a_{i,k} = a_{i} ,0 \leq i \leq n $ and $a_{n} \neq 0$ .\\
Define $P(t) = \sum_{0 \leq i \leq n} a_{i} t^{i}$ . Then roots of $P_{k}$ converge to roots of $P$ .
In particular if roots of all polynomials  $P_{k}$ are real then also roots of $P$ are real ;
if roots of all polynomials  $P_{k}$ are real nonnegative numbers then also also roots of $P$ are real nonnegative numbers .
\efac

The following corollary of Theorem 1.5 plays crucial role in our proof of Theorem 2.2 .
\cor
\begin{enumerate}
\item
Consider  a homogeneous polynomial $p(x) , x \in R^m$ of degree $n$ in $m$ real variables which is
hyperbolic in the direction $e$.
Let $x_1,x_2,x_3$ be three $e$-nonnegative vectors and $d = x_1 + x_2 + x_3$ is $e$-positive . Assume wlog
that $p(x_1 + x_2 + x_3) = 1$ . Then there exists three symmetric positive semidefinite matrices
$A,B,C$ such that $p( a_1 x_1 + a_2 x_2 + a_3 x_3) = \det( a_1 A + a_2 B + a_3 C0$ for all real $a_1 , a_2 , a_3$.
Additionally , the roots of $a_1 x_1 + a_2 x_2 + a_3 x_3$ in the direction $d$ , i.e. the roots
of the equation $p( a_1 x_1 + a_2 x_2 + a_3 x_3 - t d) = 0$ , coincide with the eigenvalues of $a_1 A + a_2 B + a_3 C$ .
\item
Theorem 2.2 is true for $e$-nonnegative tuples \\
$({\bf X})=(x_1,...x_n) , x_i \in R^m$  consisting
of at most three distinct components , i.e the cardinality of the set $\{x_1,...x_n \}$ is at most three .
\end{enumerate}
\ecor

\prf
\begin{enumerate}
\item
Consider the following homogeneous polynomial $L(b_1,b_2,b_3)= p(b_1 x_1 + b_2 x_2 + b_3(x_1 + x_2 + x_3))$  
of degree $n$ in $3$ real variables . It follows from Theorem 1.5 that there exists two real symmetric matrices
$A$ and $B$ such that $L(b_1,b_2,b_3)= \det(b_1 A +  b_2 B + b_3 I)$ . It follows that they both positive semidefinite ,
and $C = I - A - B$ is also positive semidefinite . Take a real linear combination $ z = a_1 x_1 + a_2 x_2 + a_3 x_3 $.
Then 
\begin{eqnarray*}
&p(z - t (x_1 + x_2 + x_3)) = \\
&\det( (a_1 - a_3)A + (a_2 - a_3)B + a_3 I - t I) = \\
&= \det(a_1 A + a_2 B + a_3 C - t I) .
\end{eqnarray*}
This proves that $p( a_1 x_1 + a_2 x_2 + a_3 x_3) = \det( a_1 A + a_2 B + a_3 C)$ for all real $a_1 , a_2 , a_3$ by putting $t=0$.
And it also proves the "eigenvalues " statement .
\item
Consider $e$-nonnegative tuple $({\bf X})$ consisting of $r_i$ copies of $x_i$ , $ 1 \leq i \leq 3$ ; $r_1 + r_2 + r_3 = n$ .
Assume that $d = x_1 + x_2 + x_3$ is $e$-positive (if it is not then $M_{p}({\bf X}) = 0$ by a simple argument
based on the monotonicity of $p$-mixed forms ). It follows from the polarization formula (10) that
$$
M_{p}({\bf X}) = \sum_{1 \leq i \leq k < \infty} d_i p(t_{1,i} x_1 + t_{2,i} x_2 + t_{3,i} x_3) ,
$$
and this formula is universal , i.e. holds for all homogeneous polynomial of degree $n$ , in particular for
$\det(X)$ , $X$ is $n \times n$ symmetric matrix . Therefore , using the first part of this Corollary we get that
the $p$-mixed form $M_{p}({\bf X}) = D({\bf A})$ , where the matrix tuple ${\bf A}$ consists of  $r_1$ copies of $A$ ,
$r_2$ copies of $B$ and $r_3$ copies of $C$ and $D({\bf A})$ is the mixed discriminant .
Using Rado theorem for mixed discriminants we get that $D({\bf A})> 0$ iff 
$$
Rank( \sum_{i \in S} A_i )  \geq \sum_{i \in S} r_i \mbox{ \ for all \ } S \subset \{1,2,3\} .
$$
But from the first part we get that $Rank( \sum_{i \in S} A_i )$ is equal to $p$-rank $Rank_{p}(\sum_{i \in S} x_i)$\\
of $\sum_{i \in S} x_i$
for all $S \subset \{1,2,3\}$ .
\end{enumerate}
\eprf

\pro
Consider similarly to part 4 of Fact A.1 the polynomial $p_{d}(x) =:  M_{p}(x,x,...,x,d)$ where $d$ is $e$-nonnegative
and $Rank_{p}(d) \geq 1 $ . Then $p_{d}$ is hyperbolic in any direction $z \in N_{e}(p)$ which is $e$-nonnegative and satisfies the following
inequalities :
\beqn
Rank_{p}(z) \geq n-1 ; \  Rank_{p}(z+d) = n  .
\eeqn
Also , if $y \in N_{e}(p)$ is $e$-nonnegative then also $y \in N_{z}(p_{d})$ , i.e. is $z$-nonnegative respect to the polynomial $p_{d}$.

\epro
 
\prf
Let $z \in N_{e}(p)$ be $e$-nonnegative vector satisfying (24) . Consider univariate polynomial $P(t) =  M_{p}(t z + x,t z +x,...,t z +x,d)$ .
Then $P(t) = \sum_{0 \leq i \leq n-1} a_{i} t^{i}$ and $a_{n-1} = M_{p}( z ,z ,...,z ,d)$ . It follows from Corollary  A.3 that
$a_{n-1} > 0$ . Consider now a sequence of univariate polynomials $P_{k}(t) =  M_{p}(t z_{k} + x,t z_{k} +x,...,t z_{k} +x,d_{k})$ .
Where $z_{k} ,d_{k}$ are $e$-positive and $\lim_{k \rightarrow \infty} z_{k} = z , \ \lim_{k \rightarrow \infty} d_{k} = d$ .
Then the coefficients of polynomials $P_{k}$ converge to the   coefficients of the polynomial $P$ .
It follows from part 4 of  Fact A.1 that the roots of $P_{k}$ are real . Since $a_{n-1} > 0$ hence using
Fact A.2 we get that the roots of $P$ are also real . This exactly means that the polynomial $p_{d}$ is hyperbolic in
direction $z$ . The $d$-nonnegativity statement follows from the nonnegativity part of Fact A.2 . 
\eprf

We are ready now to present our proof of Theorem 2.2 . The proof is by induction in the degree $n$ .
The main trick which we used is that to justify the induction , i.e. that if the
generalized Rado conditions hold for hyperbolic polynomial $p$ of degree $n$ then 
the generalized Rado conditions hold for some auxillary hyperbolic polynomial $p_{d}$ of degree $n-1$ ,
we need to prove  Theorem 2.2 for tuples consisting of at most three distinct components .
And this particular case follows from the classical Rado theorem via Theorem 1.5 and Corollary A.3 . 

\prf ({\bf Proof of Theorem 2.2} ) . \\

The "only if" part is simple . Indeed supposed that there exists a subset $S \subset \{1,2,...,n\}$ such that
$Rank_{p}(\sum_{i \in S} x_i) < |S| $ , i.e. using the identities (14) $M_{p}(k,k,...k,d,..,d) = 0$ ,
where $k=\sum_{i \in S} x_i$ , $d \in C_{e}(p)$ is $e$-positive and the $n$-tuple $(k,k,...k,d,..,d)$ consists of
$|S|$ copies of $k =\sum_{i \in S} x_i$ . Let $d$ be any $e$-positive positive vector such that $d - x_i$ is
$e$-nonnegative , $1 \leq i \leq n $ . Using the monotonicity of $p$-mixed forms we get that
$$
M_{p}(x_1,...,x_n) \leq  M_{p}(k,k,...k,d,..,d) = 0 .
$$
Our proof of the "if" part is by induction in the degree $n$ . 
Suppose that the generalized Rado conditions (15) hold . Then at least $Rank_{p} (x_{n}) \geq 1$ .
Consider the following homogeneous polynomial of degree $n-1$ :
$$
p_{d}(x) = M_{p}(x,x,...,x,d) , \ d = x_n .
$$
We get from Proposition A.4  the following assertion : \\
The polynomial $p_{d}(x)$ is hyperbolic in direction $z = \sum_{1 \leq i \leq n-1} x_i$ and the vectors
$x_i \in N_{z}(p_{d}) ,1 \leq i \leq n-1 $ , i.e. are $z$-nonnegative respect to the polynomial $p_{d}$.\\
Indeed , it follows from the generalized Rado conditions (15) that $Rank_{p}(z) \geq n-1$ and 
$Rank_{p}(z + d) =Rank_{p} (\sum_{1 \leq i \leq n} x_i) = n$ . \\

Next we show that the $n-1$-tuple  ${\bf Y} = (x_1,...,x_{n-1})$ satisfies the generalized Rado conditions
for $z$-hyperbolic polynomial $p_{d}$ of degree $n-1$ :
$$
Rank_{p_{d}}(\sum_{i \in S} x_i) \geq |S| \mbox{ \ for all \ } S \subset \{1,2,...,n-1 \} .
$$
Or equivalently (see formulas (14) ) , that 
\begin{eqnarray}
&M_{p}(k,..,k,z,...,z,d) > 0 ; k = \sum_{i \in S} x_i ,\\
&z =\sum_{1 \leq i \leq n-1} x_i , d = x_n , S \subset \{1,...,n-1 \} ,
\end{eqnarray}
where the $n$-tuple ${\bf T} = (k,..,k,z,...,z,d)$ consists of $|S|$ copies of $k$ , $n-1 -|S| $ copies of $z$ and one copy of $d$ .\\
It is easy to see that the generalized Rado conditions for the $n$-tuple ${\bf T}$ are implied by
the generalized Rado conditions for the original $n$-tuple ${\bf X} = (x_1,...,x_{n-1},x_n)$ .
Since the $n$-tuple $(k,..,k,z,...,z,d)$ consists of at most three distinct components hence we can
apply part 2 of Corollary A.3  . Therefore we get  that indeed \\
$M_{p}(k,..,k,z,...,z,d) > 0$ and
hence the following inequalities hold :
\beqn
Rank_{p_{d}}(\sum_{i \in S} x_i) \geq |S| \mbox{ \ for all \ } S \subset \{1,2,...,n-1 \} .
\eeqn 
Thus , by induction in the degree , we get that $p_{d}$-mixed form $M_{p_{d}}(x_1,...,x_{n-1}) > 0$ : \\
the polynomial $p_{d}$ of degree $n-1$ in $m$ real variables is $z$-hyperbolic . But\\

$M_{p_{d}}(x_1,...,x_{n-1}) = \frac{\partial^{n-1}}{\partial \alpha_1...\partial \alpha_{n-1}} p_{d}(\sum_{1 \leq i \leq n-1} \alpha_i x_i) = \\
= \frac{\partial^{n-1}}{\partial \alpha_1...\partial \alpha_{n-1}} M_{p}(\sum_{1 \leq i \leq n-1} \alpha_i x_i,...,\sum_{1 \leq i \leq n-1} 
\alpha_i x_i, x_n) = (n-1)!  M_{p} (x_1,...,x_n) .$ \\

We conclude that if Theorem 2.2 is true for $n-1$ then it is also true for $n$ , and the case "$n=1$" is trivially true .
\eprf
\rem
Consider a mixed discriminant $D({\bf A})$ , where ${\bf A} = (A_1,...,A_n)$ is a $n$-tuple of positive semidefinite $n \times n$ hermitian matrices , i.e. $A_i \succeq 0$ .
Recall that in this case $D({\bf A}) \geq 0$ ; and $D({\bf A}) > 0$ if and only if
there exists $n$ linearly independent vectors $ v_1,...,v_n$ such that $v_i \in Im(A_i) , 1 \leq i \leq n$. \\
In the proof of Theorem 2.2 we encountered the following tuple of positive semidefinite matrices : \\
${\bf A} = (A,...,A ,B,...,B,C)$ consisting of $l$ copies of $A$ , $m$ copies of $B$ and one copy of $C$ .
Moreover , this tuple is even more special . I.e. $ B-A \succeq 0$ , $Rank(B) \geq n-1$ , $Rank(A) \geq l$ ,
$rank(C) \geq 1$ ,$Rank(A+C) \geq l+1$  and $Rank(B+C) =n$ . \\
For such tuples the Rado theorem has very elementary proof , which we sketch below .\\

There are two cases . First case is when $Rank(B) = n$ , it is simple and left to the reader .\\
Second case is when $Rank(B) = n-1$ . \\
This is how in this case we can choose vectors $v_n \in Im(C) ; v_1,...,v_l \in Im(A) ; v_{l+1},...,v_{n-1} \in Im(B)$ 
in such a way that $( v_1,...,v_n )$ is a basis : first choose nonzero $v_n \in Im(C)$ which does not
belong to $Im(B)$ , second choose any $l$ linearly independent vectors $v_1,...,v_l \in Im(A$  ,
third choose any  $n-l-1$ linearly independent vectors in $Im(B) \cap L(v_1,...,v_l)^{\perp}$ .
($L(v_1,...,v_l)^{\perp}$ is a linear subspacespace orthogonal to the linear subspace $L(v_1,...,v_l)$ which
is spanned by $(v_1,...,v_l)$.) \\
\erem

\section{ Proof of Proposition 2.6 } 
\prf
Assume wlog that $q(\alpha_1,...,\alpha_n) = 1$ . It follows from the Euler's identity that 
$$
\sum_{1 \leq i \leq n}\alpha_{i } \frac{\partial}{\partial \alpha_{i} } q(\alpha_1,...,\alpha_n) = n  .
$$
Let $q(\alpha_1,...,\alpha_n) = \sum_{ (r_1,...,r_n) \in supp(q) } a_{(r_1,...,r_n)} \prod_{1 \leq i \leq n} \alpha_{i}^{r_{i}}$ . \\
Define the following nonnegative real numbers : 
$$
b_{(r_1,...,r_n)} = a_{(r_1,...,r_n) } \prod_{1 \leq i \leq n} \alpha_{i}^{r_{i}} , (r_1,...,r_n) \in supp(q) .
$$
Then $\alpha_{i } \frac{\partial}{\partial \alpha_{i} } q(\alpha_1,...,\alpha_n) =\sum_{ (r_1,...,r_n) \in supp(q) } r_{i} b_{(r_1,...,r_n)}$ .\\

Suppose that for some subset $S \subset \{1,2,...,n \} , 1 \leq |S| < n$ we have the inequality 
$\sum _{i \in S} r_i  < |S|$  for all  $(r_1,...,r_n) \in supp(q)$ . Then 
$\sum_{i \in S} \alpha_{i } \frac{\partial}{\partial \alpha_{i} } q(\alpha_1,...,\alpha_n) \leq |S|-1 $ . But the condition (17)
says that $\alpha_{i } \frac{\partial}{\partial \alpha_{i} } q(\alpha_1,...,\alpha_n) = 1 + \delta_{i}$ and
$\sum_{1 \leq i \leq n} | \delta_{i}|^{2} \leq \frac{1}{n}$ . By the Cauchy-Schwarz inequality ,
$\sum_{i \in S} | \delta_{i}| \leq \sqrt{\frac{|S|}{n}} < 1$ . Therefore ,
$$
\sum_{i \in S} \alpha_{i } \frac{\partial}{\partial \alpha_{i} } q(\alpha_1,...,\alpha_n) \geq |S| - \sum_{i \in S} | \delta_{i}| > |S|-1 .
$$
The last inequality gives a contradiction .
\eprf 

\section{ A sketch of a proof of Corollary 2.5 }
\prf
By Theorem 2.2 the conditions (1) and (2) are equivalent . (2) implies (3) for any homogeneous polynomial 
with nonnegative coefficients . \\
Let $ \alpha_i = e^{y_{i}} , 1 \leq i \leq n ; \sum_{1 \leq i \leq n} y_{i} = 0 .$ 
Consider the following convex functional
$$
f(y_1,...,y_n) = \log(q(e^{y_{1}},e^{y_{2}},...,e^{y_{n}})) .
$$
Here $q(x) , x \in R^n$ is a homogeneous polynomial of degree $n$ in $n$ real variables with nonnegative 
coefficients . Then
$$
\frac{ \alpha_{i } \frac{\partial}{\partial \alpha_{i } } q(\alpha_1,...,\alpha_n) } {q(\alpha_1,...,\alpha_n)} =
\frac{\partial}{\partial y_{i } }f(y_1,...,y_n) , 1 \leq i \leq n .
$$

Notice the condition (3) is equivalent to the following condition :
$$
\inf_{y_1+...+y_n =0} f(y_1,...,y_n) = L > -\infty .
$$ 
Consider the anti-gradient flow , i.e. the system of differential equations 
$$ 
y_{i}(t)^{\prime} = -(\frac{\partial}{\partial y_{i } }f(y_1,...,y_n)-1) , y_{i}(0) = 0 ; 1 \leq i \leq n .
$$
It is well known that in this convex case the gradient flow is defined for all $t \geq 0$ .
Using the Euler's identity we get that 
$$
\frac{d}{dt} f(y_{1}(t),...,y_{n}(t)) = -\beta(t) = :
- \sum_{1 \leq i \leq n} |\frac{ \alpha_{i } \frac{\partial}{\partial \alpha_{i } } q(\alpha_1,...,\alpha_n) } {q(\alpha_1,...,\alpha_n)} - 1|^{2}
$$
It is easy to see that , because of the convexity of $f$ , a nonnegative function $\beta(t)$ is non-increasing on $[0 , \infty )$ .\\
As $\inf_{y_1+...+y_n =0} f(y_1,...,y_n) = L > -\infty$ thus $\int_{0}^{\infty} \beta(t) dt < \infty$ .
Thus $\lim_{t \rightarrow \infty }\beta(t) = 0$ . 
This proves the implication
$(3) \rightarrow (4)$ for all homogeneous polynomials of degree $n$ in $n$ real variables with nonnegative 
coefficients  . \\
 The implication $(4) \rightarrow (5)$ is obvious .  The implication $(5) \rightarrow (6)$ for general 
homogeneous polynomials of degree $n$ in $n$ real variables with nonnegative coefficients is Proposition 2.6 . \\
Finally , the implication  $(6) \rightarrow (2)$ follows fairly directly from Theorem 2.2 . 
\eprf

\section {Lower bounds on the number of oracle calls for the exact computation of $\frac{\partial^n}{\partial x_1... \partial x_n} p(x_1,...,x_n)$ }
\dfn
Call a set $\{X_1,...,X_m\} , X_i \in C^{n}$ $\epsilon$-universal if there exist complex numbers $c_1,...,c_m$ such that for any homogeneous polynomial $p(.)$ of degree $n$ in $n$ complex variables 
the following inequality holds
\begin{eqnarray}
&|\frac{\partial^n}{\partial x_1... \partial x_n} p(x_1,...,x_n) -\sum_{1 \leq i \leq m } c_{i} p(X_{i})| \\ \nonumber
& \leq \epsilon \max_{(r_1,...,r_n) \in I_{n,n}} |a_{r_1,...,r_k}| ,
\end{eqnarray}
where $a_{r_1,...,r_n} ,(r_1,...,r_n) \in I_{n,n} $ are the coefficients of the polynomial $p(.)$ .
\edfn
\lem
If the set $\{X_1,...,X_m\} , X_i \in C^{n}$ is $0$-universal then 
\beqn
m \geq \frac{n!}{[\frac{n}{2}]! (n-[\frac{n}{2}])!} \approx \frac{2^{n}}{\sqrt{n}}
\eeqn
If the set $\{X_1,...,X_m\} , X_i \in C^{n}$ is $\epsilon$-universal then 
\beqn
m \geq \min([ \frac{1}{\epsilon}] ,\frac{n!}{[\frac{n}{2}]! (n-[\frac{n}{2}])!} ) 
\eeqn
\elem
\prf
Define a monomial $M_{r_1,...,r_n}(x_1,...,x_n) = x_{1}^{r_1} x_{2}^{r_2}...x_{n}^{r_n}$ .
As $\{X_1,...,X_m\}$ is universal thus the exists complex numbers $(c_1,...,c_m)$ , which are wlog are all nonzero , such
that 
$$
\sum_{1 \leq i \leq m} c_i M_{r_1,...,r_n}(c_{i}^{\frac{1}{n}} X_{i}) = 0
$$ 
if $(r_1,...,r_n) \in I(n,n) ,(r_1,...,r_n) \neq (1,1,...,1) $ ; \\
and $\sum_{1 \leq i \leq m} M_{1,1,...,1}(c_{i}^{\frac{1}{n}} X_{i}) = 1$ ; define $Y_i = c_{i}^{\frac{1}{n}} X_i$ (here $c_{i}^{\frac{1}{n}}$ is one of the $n$th complex roots  of $c_i$ ).\\   
Let $Half = \{(r_1,...,r_n) : r_i \in \{0,1\} , \sum_{1 \leq i \leq n} r_i = [\frac{n}{2}]$ . Notice that the cardinality
$|Half| = \frac{n!}{[\frac{n}{2}]! (n-[\frac{n}{2}])!} = : K$ . \\ 
Define the following two $K \times m$ complex matrices :
$$
W((r_1,...,r_n) , j) = M_{r_1,...,r_n} (Y_{j}) ,
$$
 
\begin{eqnarray*}
&V((r_1,...,r_n) , j) = M_{1-r_1,...,1-r_n} (Y_{j}) : \\
&(r_1,...,r_n) \in Half , 1 \leq j \leq m  .
\end{eqnarray*}
Clearly , $Rank(W)=Rank(V) \leq m$ . On the other hand the $0$-universality condition implies the matrix identity
\beqn
W V^{T} = I 
\eeqn
Therefore $m \geq Rank(W) \geq |Half| = \frac{n!}{[\frac{n}{2}]! (n-[\frac{n}{2}])!} $ . \\
If the set $\{X_1,...,X_m\}$ is $\epsilon$-universal and $d \epsilon < 1 , d \in N$ then $Rank(W V^{T}) \geq d$ .
This proves (30) .
\eprf
\rem
The identity (7) is a particular case of a slightly more general one :
\beqn
\frac{\partial^n}{\partial x_1... \partial x_N} p(x_1,...,x_n)= E ( p(z_{1},z_{2},...,z_{n}) \prod_{1 \leq i \leq n } \overline{z_{i}} ),
\eeqn
where $(z_{1},z_{2},...,z_{n})$ are independent complex random variables such that $E(z_{i})=0$ and $E(z_{i}\overline{z_{i}}) = 1$ for
all $1 \leq i \leq n$ . The identity is easily proved by checking it for all monomials $M_{r_1,...,r_n} , (r_1,...,r_n) \in I(n,n)$ .
If $p(.)$ is a multilinear polynomial , i.e.
$$
p(x_1,...,x_n) = \prod_{1 \leq i \leq n} (\sum_{1 \leq j \leq n} A(i,j) x_j ) ,
$$
then $\frac{\partial^n}{\partial x_1... \partial x_N} p(x_1,...,x_n) = Per(A)$ , where $per(A)$ is the permanent of the matrix $A$ .
Clearly , lower bound $m \geq \frac{n!}{[\frac{n}{2}]! (n-[\frac{n}{2}])!}$ also holds for multilinear polynomials and even for
powers $(\sum_{1 \leq i \leq n} a_i x_i )^{n}$ . It is very likely that the actual lower bound is $2^{n-1}$ and that
it does exist somewhere in the geometrical designs literature . In the case of permanents , the formula (7) is essentially the 
Ryser's formula \cite{minc}; and Lemma D.2 says that , in some sense , it is an optimal formula for computing permanents . \\
Another equivalent formulation of Lemma D.2 is the following statement : \\
Let a set of complex vectors 
$$
S = \{X_l =(x_{l,1},...,x_{l,n}) \in C^{n} : 1 \leq l \leq \frac{(2n-1)!}{(n-1)! n!} \}
$$ 
be a Haar set for the monomials $M_{r_1,...,r_n} : (r_1,...,r_n) \in I_{n,n}$. \\
I.e. the square matrix $\{M_{r_1,...,r_n}(X_i) : X_i \in S ; (r_1,...,r_n) \in I_{n,n} \}$ is nonsigular . \\ 
If 
$$
\prod_{1 \leq i \leq n} x_{l,i} = \sum_{1 \leq k \leq m} c_i (\sum_{1 \leq i \leq n} Y(k,i) x_i)^{n}
$$ 
for all $1 \leq l \leq \frac{(2n-1)!}{(n-1)! n!}$ and
some complex numbers $\{ c_k ;Y(k,i) :1 \leq k \leq m , 1 \leq i \leq n \} $ then 
$$
m \geq \frac{n!}{[\frac{n}{2}]! (n-[\frac{n}{2}])!} \approx \frac{2^{n}}{\sqrt{n}} .
$$
\erem

\end{document}